# Fighting with deterministic disturbances

V. N. Tibabishev

Asvt51@narod.ru

**Consider the problem of interference mitigation in the identification of the dynamics of multidimensional control systems in the class of linear stationary models for single realizations of the observed signals. Instead unverifiable concepts uncorrelated processes, the notion of sets of components of the signal measured on a semiring. Properties of signals are defined for systems of sets of linearly dependent and linearly independent measured signals. Frequency method is found to deal with noise on the set of deterministic functions. Example is considered to identify the dynamic characteristics of the aircraft on the data obtained in the regime of one automatic landing.**

1. Introduction.

Known methods of statistical dynamics control systems cannot always be used, for example, when solving the problem of identification of dynamic characteristics of the aircraft. Mathematical model of the dynamics of the aircraft are difficult to obtain analytical techniques due to elastic deformation the airframe. Therefore, the dynamic characteristics are found by solving the problem of identification from experimental data obtained, for example, in one mode, the automatic landing class of IL-96. The aircraft is a multidimensional subject of management. Dynamics of an aircraft was necessary to determine only two control channels. The first control channel is specified by pitch - $x_1$- pitch of the aircraft – $y$. On the pitch the aircraft affect the position of the flaps $x_3$ and the position angle of the stabilizer $x_4$. In one mode, automatic landing were received simultaneous recording inputs $x_i$, $i = 1,2,3,4$ and the output signal from 274 reference in 0.5 sec.

A well-known method of solving this problem [1, p.22] is reduced to solving a system of integral equations of the first kind. The system of integral equations contains all the autocorrelation functions of the input signals, all the cross-correlation function between the inputs and all the cross-correlation function between all input and output signals. In the process of automatic landing is performed changing the angle of the stabilizer and additional release flaps. For this reason, the actual pitch of the aircraft is a no stationary random process. It is known [2, pp. 458] that the non-stationary random processes do not possess the ergodic property. Therefore, the correlation function obtained by properly averaging over time of individual implementations possible. Correlation functions obtained by averaging over the set of realizations of the same is not possible, as impossible to make a set of landing an aircraft under the same conditions, such as weather. In addition, the system should have four independent outputs (exits). In this case there is only one way out. The task cannot be solved by well-known for these reasons.

It is known [2, pp. 374] that if the realization of a random process is defined as a single implementation, it is an ordinary (nonrandom) function. Thus, the problem of solving the problem of identification of dynamic characteristics of control objects is not on the set of random functions, and on the set of single deterministic functions.

The aircraft is a nonlinear and no stationary object of control [3, pp. 420]. However, in a small interval of observation with an almost constant velocity flight on final approach the aircraft can be assumed linear and stationary objects of management. Therefore, we will solve the problem of identification of dynamic characteristics in the first approximation in a class of linear stationary models.

The task of identifying the dynamic characteristics precedes the solution of the problem of suppression of interference. The problem of noise reduction can be solved if the desired signals and noise differ in their properties. Come to the necessity of constructing mathematical models of signals and noises in the class of nonrandom functions.

## 2. Measurable processes and signals

Management systems may profit by one-dimensional and multi-dimensional, unsteady and nonlinear. In the simplest case, the dynamics of one-dimensional control object in a first approximation is in a class of linear stationary models, which describe the differentiation operator with constant coefficients, for example $n$ - the second order $Dy = x$. When solving practical problems rather than the exact original data can be obtained only measured input data. Under the exact data we mean the components of the measured signals in the form of functions related to operator equations. In reality, accurate baseline data observed in additive and multiplicative noise, which are not linked by an operator equation. Noises at the input do not generate involuntary movement, and noise output is not part of the forced motion in the channel management. Therefore, interference and accurate data are of different nature.

All that is described as a function of time is called a process. Process, converted into electrical form (analog or digital) is convenient for further processing, is called a signal. Management systems may have one or more control channels. Channel Management - the device displays a set of precise components of the input signals (a) set the exact components of the output signals. A mathematical model of the control channel is determined by the operator. Operators can be linear and nonlinear stationary and no stationary. Multidimensional control system is divided into systems with linearly independent and linearly dependent on input actions. Initially, we consider the multidimensional system with linearly independent input actions. In considering the general properties of the input $x(t)$ or output $y(t)$ processes or signals will be denoted by their generalized symbol.

Further, we consider cases where the observed signal is given only a single implementation, taken from the set of signals. In the real world instead of the exact components of the signals $z$, coupled operator equations that can be measured in (recorded) only approximate the signals $\tilde{z}$, which are distorted by additive and multiplicative noise $\tilde{z} = \vartheta z + n$, where $\vartheta$ and $n$ - respectively the multiplicative and additive noise. The observed signal is reduced to

$$\tilde{z} = \hat{\vartheta} z + p + n, \qquad (1)$$

where $\hat{\vartheta}$ - a number equal to the average value of the multiplicative noise $\vartheta$,

$p$ -reduced additive noise generated by multiplicative noise. If the multiplicative noise is absent, we assume that $\hat{\vartheta}$ =1, and $p$=0. If the approximate observed signals of absence related to an operator equation, z=0.

Given to (1) the processes or signals of any nature will be called measurable processes or signals. In many cases, the exact input action is a result of the interaction of material bodies, for example, translational motion steering column of the IL-62, occurs under the force action pilots. The basis of the mechanics of interacting material points is the Newton equation of motion [5, pp. 19], according to which body mass $m$, the second derivative $\ddot{x}$ of the function, for example, describing the translational motion $x$ is connected with the power $\gamma$ of ordinary inhomogeneous differential equation $m\ddot{x} = \gamma$. In flight, the movement $x$ may vary not only by the will of the pilot, as a result of, for example, random wind gusts. In addition to the wire system of air traffic steering column sensor is converted into electrical signals. When converting processes into electrical signals and their transmission in addition there are additive and multiplicative noises.

The differential equation of Newton is a linear differential equation. Therefore, under the symbol x can be any component of the measured signal. For each component of the measured signal (1), for example, n condition for the existence and uniqueness of solutions of differential equations of motion $m\ddot{n} = \gamma_n$ performed. It follows that $\gamma$ and the second derivative $\ddot{x}$ is a bounded continuous function [6, pp. 150]. Obviously, the first derivative $\dot{x}$ of the function $x$ itself is also a continuous function.

For the harmonic effects $\gamma_k = A_k exp(j\omega_k t)$ with arbitrary frequency $\omega_k$ and bounded amplitude $A_k$ there is a solution of the inhomogeneous differential equation of Newton $x(t) = -A_k \exp(j\omega_k t)/\omega_k^2$ on the whole real line $t \in (-\infty, +\infty)$. Suppose the right side of the differential equation of Newton is an optionally convergent trigonometric series. It is known [6, pp. 234] that this differential equation will be presented to the convergent trigonometric series. It is known [7. 444] that if $\gamma$ is a continuous function on the whole line, it just seems to be convergent trigonometric series. Further, we believe that elements of the measured processes $z, p$ and $n$ belong to the set of bounded continuous functions can be represented by convergent trigonometric series. On the set of functions square, for example, the input action is not integrable on the entire line of Riemann and Lebesgue and Stieltjes integrable function on $g(t) = t/2T$

$$F = \lim_{T\to\infty} \int_{-T}^{+T} x_k^2(t)dg(t) = \lim_{T\to\infty} \frac{1}{2T}\int_{-T}^{+T} x_k^2(t)dt \equiv M\{x_k^2(t)\} < \infty. \quad (2)$$

These properties are, for example, functions that belong to the Hilbert space of almost periodic functions $B_2(-\infty, +\infty)$ in the sense of Besicovitch [12, c. 222]. Further, we assume that $z, p, n \in B_2(-\infty, +\infty)$.

Proper motion in the channel management is a solution of the homogeneous differential equation $m\ddot{x} = 0$. Such a solution satisfies the differential equation $x_1(t) = ct + b$, where c and b - arbitrary bounded constants. Not difficult to see that an arbitrary value of $c \neq 0$ even if $b = 0$ proper motion system does not satisfy (2). Only when $c \equiv 0$ and arbitrary $b = 0$ natural and forced solution of differential equations of Newton satisfies the limited power of the signal or (2). If the arbitrary constant $b$ equal to zero, we find that the general solution of differential equations of Newton coincides with the partial solution of inhomogeneous differential equation.

All the components of the system of sets of measured input and output of multidimensional control systems can be divided into sub-system of sets of linearly independent components of the measured signal and noise $\mathbb{R}$, and the subsystem is linearly dependent components of the exact input and output signals $\mho$, coupled operator equations.

## 3. Subsystem is linearly independent sets of processes

Let the control system has $d$ inputs and $h$ outputs. On the sets of measured signals can distinguish several sub-linearly independent sets of signals. For example, the input $k$ measured signal $\tilde{x}_k = \hat{\theta} x_k + p_k + n_k$ belongs to the system sets $\tilde{x}_k \in \mathbb{R}_{xk} = X_k \cup P_k \cup N_k$, where $\hat{\theta} x_k \in X_k$ a set of precise components of the signals that generate the forced movement of the control channels. $P_k$ - the set of reduced additive noise generated by the multiplicative interference $\theta$, where $\hat{\theta}$ - the average value. $N_k$ -set of additive noise, $k$ -dead entrance. Another system of sets of linearly independent elements is the union of all input signals $\mathbb{R}_x = \uplus_{k=1}^{k=d} \mathbb{R}_{xk}$. Similarly, we obtain a system of linearly independent sets of output signals $\mathbb{R}_y = \uplus_{k=1}^{k=h} \mathbb{R}_{yk}$, where $\tilde{y}_k \in \mathbb{R}_{yk} = Y_k \cup P_{yk} \cup N_{yk}$, $\tilde{y}_k = \hat{\vartheta} y_k + p_{yk} + n_{yk}$. $Y_k$ - a set of accurate output signals. $P_{yk}$ - the set of reduced additive noise generated by the multiplicative interference $\vartheta$, where $\hat{\vartheta}$ - the average value. $N_{yk}$ - set of additive noise, $k$ -died out. This system of linearly independent sets of elements observed at the output can be enhanced by add to it the system of sets of linearly independent components of interference, which distort the input signals.

Each subset of elements, for example, $\mathbb{R}_{xk}$ conceived as a whole in the sense that each element is a subset of the linear hull spanned by a particular orthonormal basis of harmonic $f_{kl} \in G_k$, $l = 1,2,3,...$, $f_{kl} = \exp(j\omega_{kl} t)$. In addition, each $i$ - th element of a valid, for example, the exact component of the input signal at the input $k$ is represented as

$$x_{ki} = \sum_{l=1}^{l=c}(a_{kil} \cos(\omega_{kl} t) + b_{kil} \sin(\omega_{kl} t)). \qquad (3)$$

If the coefficients $a_{kil}$ and $b_{kil}$ are uncorrelated with each other with random numbers, depending on the number of realization of $i$, then the elements $x_{ki}$ are stationary in the wide sense stochastic processes, which, when averaged over multiple realizations generate a correlation function of the general form $R_k(\tau) = \sum_l A_{kl} \cos(\omega_{kl} t)$ [8, 336].

Usually in scientific articles suggest that the additive noise distorted signals belong to the same functional space. In the functional spaces for the elements of the sets introduced the concept of norm. However, the observed signals are the union of different properties of sets of precise components and noise. In general, each component of the exact signal closely observed in the presence of two additive noises belonging to different subsets. For the elements of the sets introduced the measure $\mu$ (V) - numerical characteristic satisfying:
1) the domain of $\mu$ (V) is a semiring of sets,
2) the values of $\mu$ (V) are valid and non-negative,
3) $\mu$ (V) is additive, i.e., for any finite decomposition $V = V_1 \cup ... \cup V_k$ into disjoint sets, the equality $\mu(V) = \sum_{l=1}^{c} \mu(V_l)$ [9, sec. 249].

It is known [9, sec. 235] that a measure can be taken, for example, the integral of a nonnegative function. Take as a measure of the sets of measured signals of the functional (2), the domain which should be a semiring of sets. Subsystem sets of functions $\mathbb{R}_d$ a semiring if it is a union $\mathbb{R}_x = \uplus_{k=1}^{k=d} \mathbb{R}_{xk}$ disjoint sets satisfying $\mathbb{R}_{xk} \cap \mathbb{R}_{xm} = \emptyset$ for all $k \neq m$, where $\emptyset$ - empty symbol set.

Let $x_{ki} \in \mathbb{R}_{xk}$ arbitrarily selected item from the subset $\mathbb{R}_{xk}$, $1 \leq k \leq d$, $i = 1,2,3...$, $\tilde{x}_{ki} = \hat{\vartheta} x_{xki} + p_{xki} + n_{xki}$. Functional (2), defined for an arbitrarily selected elements of the union of sets $\mathbb{R}_x = \cup_{k=1}^{k=d} \mathbb{R}_{xk}$, $F(V) = M\{\sum_{k=1}^{k=d} \tilde{x}_{ki}\}^2$ will be additive numerical characteristic, if all the scalar products for the randomly selected $(\tilde{x}_{ki}, \bar{\tilde{x}}_{mi}) \equiv M\{\tilde{x}_{ki} \cdot \bar{\tilde{x}}_{mi}\} = 0$ for all $k \neq m$, where $\tilde{x}_{ki}$ and $\bar{\bar{x}}_{mi}$ - complex conjugate functions. Known [9, p.136] that if the elements are linearly independent, then they are orthogonal. In determining the functional (2) all scalar products are found in pairs of linearly independent and hence orthogonal functions. This implies that the elements are linearly independent systems of sets in Hilbert space of almost periodic functions are a semicircle, and the functional (2) defines a measure on it.

Check the condition $\mathbb{R}_{xk} \cap \mathbb{R}_{xm} = \emptyset$ on infinite set of processes is almost impossible. Let $V_k, V_m$ are metric spaces and contain $n$ elements. For each element $x_i \in V_k$ and every item $y_l \in V_m$ can define $n^2$ values of the distances $\rho(x_i, y_l)$. If for all , $l = 1,2,3, \ldots n$, the condition $\rho(x_i, y_l) \neq 0$, then the suppression $V_k \cap V_m$ is empty. Such an algorithm is impossible to implement in practice, because each subset $\mathbb{R}_{xk}$ contains an infinite countable set of elements $n = \infty$. In this connection there is a need to develop realizable in practice, the method for determining the semiring for a system of infinite sets of components of the observed signals.

Theorem 1. *Given a system of infinite sets of functions* $\mathbb{R} = \cup_{k=1}^{k=q} V_k$, *where q - the number of components of the measured signal is the linear span* $V_k = \{x_{ik} \in V_k : x_{ik} = \sum_p a_{pki} e_{pk}\}$, *spanned by the corresponding orthonormal bases* $e_{pk} = \exp(j\omega_{pk} t)$, $p, i = 1,2,3, \ldots, p, i = 1,2,3, \ldots$ *in the Hilbert space of almost periodic functions. If the system of sets of functions* $\mathbb{R}$ *is a semiring, then the numerical sets of angular frequency* $\omega_{pk} \in \Omega_k$ *harmonic bases, forming a system of infinite sets of functions of* $\mathbb{R}$, *is a semicircle of disjoint sets of numerical angular frequency of the harmonic components of orthonormal bases* $\Omega_k \cap \Omega_l = \emptyset$ *for all* $k \neq l$.

It is proof. Since $\mathbb{R}$ is a semiring, then the functional (2) will define the additive measure if for all $x_{ki} \in V_k, x_{pl} \in V_p$, $p, i, l = 1,2,3, \ldots, k = 1,2,3, \ldots q$ inner products are zero for all k ≠ p. Scalar product is zero for all $a_{ki}\bar{b}_{pl}$ for arbitrarily chosen i and , if for all k ≠ p be the condition $\omega_k - \omega_p \neq 0$

$$M\{(\sum_i^c a_{ki} \exp(j\omega_k t) \sum_l^c \bar{b}_{pl} \exp(-j\omega_p t))\} = M\{(\sum_i^c \sum_l^c a_{ki} \bar{b}_{pl} \exp j(\omega_k - \omega_p) t)\} = 0. \quad (4)$$

This implies that the system of numerical sets of angular frequency $\mathbb{Q} = \cup_{k=1}^{k=d} \mathbb{Q}_k$ is a semiring where $\mathbb{Q}_k = \cup_{kc=1}^{c=q} \Omega_{kc}$, where $k$ - number of input, $q$ - the number of independent components of the measured signal.

It is consequence from theorem 1. It is known [10, pp. 186] that the Hilbert space of almost periodic functions defined by a completely continuous symmetric operator normal form of convolution type, displaying $B_2(-\infty, +\infty)$ into itself. Known [11, pp. 214] that if the operator is compact and normal form, it generates only a finite number of basic functions and non-zero eigenvalue. Therefore, in terms of elements (3) and the scalar products (4) the value of $c < \infty$. This implies that the system $\mathbb{Q} = \cup_{k=1}^{k=d} \Omega_k$ numerical sets of values of the angular frequency of harmonic basis functions of the independent components of input signals is the union of pair wise disjoint finite number of sets.

It is corollary 2 of Theorem 1. It is known [8, with. 340], what the set of almost periodic functions casual stationary process is represented by a trigonometrically series (3) with material numbers $\omega_k$. In the theory of random processes domain of (3) is only a set of almost periodic functions. Components of the measured signal as represented by series (3), but their domain is not a separate set of components and the system of sets, which should be a semicircle. In this case, the frequencies $\omega_k$ are choice the additional restrictions. Periodic processes are represented by harmonic series (3) with multiple frequencies as well belong to the set of almost periodic functions. For example, suppose that a set of frequencies includes frequencies $\Omega_a$ contains frequencies $\omega_k = k\omega_1$, $k = 1,2,3,...$, and set $\Omega_b$ contains frequencies $\omega_d = d\omega_2$, $d = 1,2,3,...$ Let $\omega_2 = \beta\omega_1$, where $\beta < 1$ - an arbitrary rational number, for example, $\beta = n/m$. When $d = m$ and $k = n$ suppression sets will not be empty. This implies that the system of linearly independent components of the processes observed cannot contain more than one subsystem, containing periodic components, which are represented by harmonic series with multiple frequencies. The condition of independence components of the measured signals will be executed always; if the system sets the frequency of random harmonic components will contain only the arbitrary incommensurate frequencies.

4. Subsystem of linearly independent sets of processes

Synchronously with the input signals the output signals are observed for all h outputs. Instead, accurate output signals $y_r$ can be observed only by their approximate values $\tilde{y}_r = \theta_r y_r + \vartheta_r y$, where $\vartheta_r$ - multiplicative noise, $\vartheta_r$ - additive noise, r - number of the exit. Given the known property of multiplicative noise [4, pp. 139], we obtain an expression for the output of the measured signal $\tilde{y}_r = \hat{\theta}_r y_r + \beta_r + \vartheta_r$, where $\hat{\theta}_r$ − mean multiplicative noise $\theta_r$, $\beta_r$ - reduced additive noise generated by the multiplicative interference $\theta_r$.

It is known [10, c. 186] that the Hilbert space of almost periodic functions defined by a completely continuous operator of convolution type, normal shape, generated by the exact component of the input signal $x(t)$, which displays $B_2(-\infty, +\infty)$ into itself. In the multidimensional control systems in general, each output additively associated with each input through a particular control channel. In this case, the exact components of the observed signals synchronously $x_d$ all d inputs and one output r are related by an operator equation

$$y_r(t) = \sum_{c=1}^{c=d} M\{x_c(t - \tau) k_{rc}(\tau)\} = \sum_{c=1}^{c=d} A_c k_{rc}, \qquad (5)$$

where $k_{rc}$ - the weight function of the control channel between the output and input r c.

Instead of the exact component $y_r$ output signal for each output signal is observed r observed signal $\tilde{y}_r$, belonging to a system of sets of linearly independent components of the output signals $\tilde{y}_r \in \mathbb{Z}_r = Y_r \cup E_r \cup T_r$, where the set of elements output signals $\hat{\theta}_r y_r \in Y_r$, $\beta_r \in E_r$, $\vartheta_r \in T_r$. The union of these sets of systems on all outputs leads to a system of sets of linearly independent components of the output signals $\mathbb{Z} = \cup_{r=1}^{r=h} \mathbb{Z}_r$. If the systems are sets of input and output signals of R and Z are semi rings, then their union loses this property.

Theorem 2. *Suppose that in certain functional spaces, the set of weight functions $k \in K$ is displayed on the set Y by a linear stationary operator of convolution type $y = Ak$, $y \in Y$, where the kernel of A is generated by the input signal $x \in X$. If the operator equation $y = Ak$ is defined in the Hilbert space of almost periodic functions in the sense of Besicovitch, the image of a completely continuous operator of convolution type is the linear span [10, pp. 50], which stretched over a finite orthonormal basis of the kernel.*

It is proof. It is known [12, pp. 224] that every almost periodic function of Besicovitch is convergent trigonometric series, for example, $x = \sum_{c=1}^{c=s} a_c b_c$, where $b_c = \exp(j\omega_c t)$. Using Fourier transform we can obtain that the operator of convolution type is represented as $y = \sum_{c=1}^{c=s} \lambda_c \alpha_c b_c$, where $\lambda_c$ - Fourier exponents of the weight function $k \in K$ define the projection of the weight function for the orthonormal function $b_c$. It follows that the image of the operator of convolution type is spanned spanned by the orthonormal basis of the kernel $b_c \in G$. Known [11, pp. 214], which are completely continuous operators generate a finite number of eigenvalues and eigenvectors. Therefore, the representations of the kernel, image and Fourier transform $< \infty$.
Theorem 2 proved.

It is consequence of Theorem 2. In the Hilbert space of almost periodic functions of the image of $y \in Y$ and a function that generates the kernel of $\in X$, a linear stationary operator of convolution type are linearly dependent elements. The linear dependence of the image $y \in Y$ and a function that generates the kernel of $x \in X$, a linear stationary operator of convolution type follows from the fact that equality $\sum_{c=1}^{c=s} p_c a_c b_c + \sum_{c=1}^{c=s} \lambda_c a_c b_c = 0$ holds for all nonzero $a_c, b_c$, $p_c$ and $\lambda_c$, when $p_c = -\lambda_c$. It is known [10, p. 185] that the Hilbert space of almost periodic functions $B_2(-\infty, +\infty)$ is an uncountable set of functions $\exp(j\omega t)$, and each vector in this space is only a countable set of nonzero eigenvectors $b_c = \exp(j\omega_c t)$. Theorem 2 implies that the finite numeric sets of frequencies $\omega_c \in \Omega_x$ and $\omega_c \in \Omega_y$ coincide $\Omega_x = \Omega_y$. It follows that ate a system of sets contains at least one image and inverse image of a linear operator, then the system sets is not a semicircle.

### 5. Subsystem of linearly dependent input of processes

Consider the method of dealing with noise in the special case when all inputs $d$ control system there are linearly independent input effects. In such cases, the system sets constitute the input signals $\mathbb{R} = \uplus_{k=1}^{k=d} \mathbb{R}_k$ is a union of disjoint sets of components of the signals. By Theorem 1, the system sets the values of the angular frequency of harmonic functions $\mathbb{Q} = \uplus \Omega_r$ is a union of disjoint finite number of sets. Interconnection of the sets $\mathbb{Z} = \cup_{r=1}^{r=h} \mathbb{Z}_r$, observed on all outputs $h$, is just a semicircle.

*Lemma. Let the system sets the input $\mathbb{R} = \uplus_{k=1}^{k=d} \mathbb{R}_k$ and output $\mathbb{Z} = \cup_{r=1}^{r=h} \mathbb{Z}_r$ the measured signals in a multidimensional control system containing $d$ inputs and $h$ outputs are semirings. If between the number of input $i$ and output number $l$ there exists a linear stationary control channel, represented by a convolution operator, defined in the Hilbert space of almost periodic functions in the sense of Besicovitch, then the intersection number of sets of angular frequency components of orthogonal bases, the observed signals at the input $i$ and output $l$ contains a finite subset of the numeric values of the angular frequency accurate component of the input signal observed at the inlet $i$.*

It is proof. Semiring of sets $\mathbb{R} = \uplus_{k=1}^{k=d} \mathbb{R}_k$ semiring is the union of systems of sets, $\mathbb{R}_k$, $k = 1,2,3,.,i,..d$. In this semiring is a semiring $\mathbb{R}_i$ components $\mathbb{R}_i = X_i \cup P_i \cup N_i$, where the infinite set of exact components of the signals $\hat{\theta}x \in X_i$, given the additive noise $p \in P_i$, additive noise $n \in N_i$. According to Theorem 1 semicircle of infinite-dimensional sets of elements $\mathbb{R}_i$ generated semiring of finite sets of angular frequency harmonic bases

$$\widetilde{\Omega}_{xi} = \Omega_{xi} \cup \Omega_{pi} \cup \Omega_{ni}. \tag{6}$$

By (5) at the output signal is observed $l$ equal to the sum of signals $y_l = \sum_{c=1}^{c=d} A_c k_{cl} = \sum_{c=1}^{c=d} y_{cl}$, where $y_{cl} \in \tilde{Y}_{cl} = \cup_{c=1}^{c=d} Y_{cl} \cup E_l \cup T_l$, which contains a semicircle $\tilde{Y}_{il} = Y_{il} \cup E_l \cup T_l \subset \tilde{Y}_{cl}$. This infinite system of sets of functions generated by a finite number system, sets the angular frequency of harmonic bases $\widetilde{\Omega}_{yil} = \Omega_{yil} \cup \Omega_{El} \cup \Omega_{Tl}$, where the finite number sets the angular frequency, respectively, $\Omega_{yil}$- exact output signal generated by the operator $A_i k_{il}$, $\Omega_{El}$- reduced additive noise, $\Omega_{Tl}$- additive noise. By Theorem 2 $\Omega_{yil} = \Omega_{xi}$. In this case, $\widetilde{\Omega}_{yil} = \Omega_{xi} \cup \Omega_{El} \cup \Omega_{Tl}$. Lemma proved.

It is consequence of the lemma. Lemma determines the order of solving the problem of precise allocation of components of the signals observed, for example, the input $i$ and output $l$ in the presence of multiplicative and additive noise. Initially determined by the numeric arrays frequency $\widetilde{\Omega}_{yil}$ and $\widetilde{\Omega}_{xi}$, using, for example, a well-known ideal analyzer [13, pp. 113] to determine the current amplitude spectrum [13, pp. 112]. Arrays frequencies are determined, for example, local maxima of the current amplitude spectra. Then there is the intersection of the found arrays of frequencies $U_{il} = \widetilde{\Omega}_{yil} \cap \widetilde{\Omega}_{xi}$. Isolation of accurate projections constitute the input and output signals observed on the background of multiplicative and additive noise is performed using the definitions of indicators of the Fourier series only at frequencies $\omega_k \in U_{il}$.

### 6. Subsystem of linearly dependent input of processes

When solving practical problems, as a rule, all or part of the input signals are linearly dependent. In such cases, the property of the lemma is not satisfied and the above described algorithm for solving applied problems is not applicable. The degree of linear relationship between input actions can be estimated by the cosine of the angle $\cos(\varphi_{il})$ [9, p. 135] between the exact components of the observed signals (prototypes) for inputs $i$ and $l$. It highlighted the following cases. If a selected pair of inputs observed transforms signals are linearly independent, then $\cos(\varphi_{il}) = 0$. If the inputs, for example, $i$ and $d$ transforms the observed signals are identical up to sign, then $|\cos(\varphi_{id})| = 1$. In this case, the dimension of the control system of inputs is too high. One input $i$ or $d$ should be excluded from consideration. After exclusion of redundant inputs, we find that for systems with linearly dependent input actions for all inputs, the condition $|\cos(\varphi_{il})| < 1$.

Theorem 3. *Let the control system has $d$ inputs. Each input $\tilde{x}_i$, $1 \leq i \leq d$ is an element of a semiring of sets $\tilde{x}_i \in \widetilde{\mathbb{R}}_i$ components of a signal containing communication features (extra preimages) $f_{il}$, $\tilde{x}_i = \hat{\vartheta}_i x_i + p_i + n_i + \sum_{l=1}^{l=d} f_{il}$ for all $l \neq i$ with a measure generated by the scalar product in Hilbert space of almost periodic functions. If at all entrances, or some part, the linear dependence of signal transforms $\hat{\vartheta}_i x_i + \sum_{l=1}^{l=d} f_{il}$, satisfying $|\cos(\varphi_{il})| < 1$ for all $i \neq l$, then from the original system of linearly dependent input signals $\widetilde{\mathbb{R}} = \cup_{i=1}^{i=d} \widetilde{\mathbb{R}}_i$, containing communication features can be a subsystem of linearly independent input influences $\hat{x}_i = \hat{\vartheta}_i x_i + p_i + n_i \in \mathbb{R} \subset \widetilde{\mathbb{R}}$, containing no communication function.*

It is proof. That would meet the conditions for prototypes of input signals $|\cos(\varphi_{il})| < 1$ for all $l \neq i$ represent the measured input signal in the form $\tilde{x}_i = \hat{\vartheta}_i x_i + p_i + n_i + \sum_{l=1}^{l=d} f_{il}$ for all $l \neq i$. This representation differs from the representation (1) the presence of additional prototypes as a liaison between the input channels $f_{il}$. By hypothesis, each signal $\tilde{x}_i$ is an element of a semiring of sets $\widetilde{\mathbb{R}}_i$ components of the measured signals on the selected input. Therefore, all possible inner products between the components of the signal $\tilde{x}_i$, for example, components of reimages $(\hat{\vartheta}_i x_i, \bar{f}_{il}) \equiv 0$ identically equal to zero. This condition is fulfilled for all $1 \leq i, l \leq d$, $l \neq i$. Union $\widetilde{\mathbb{R}} = \cup_{i=1}^{i=d} \widetilde{\mathbb{R}}_i$ is not a semicircle, as suppression systems, sets of input signals contain non-empty sets of subsystem coupling functions, such as $\widetilde{\mathbb{R}}_i \cap \widetilde{\mathbb{R}}_l = \sum_{l=1}^{l=d} f_{il}$, $l \neq i$. Let us find easy to implement in practice the symmetric difference [9, p. 13], two systems of sets $\widetilde{R}_i \Delta \widetilde{R}_l = \widetilde{\mathbb{R}}_i \cup \widetilde{\mathbb{R}}_l \setminus \widetilde{\mathbb{R}}_i \cap \widetilde{\mathbb{R}}_l = \mathbb{R}_i \cup \mathbb{R}_l$. With the help of the available operations, define the union of all symmetric differences for all $i \neq k$ $\cup_{i=1}^{i=d} (\cup_{k=1}^{k=d} \widetilde{\mathbb{R}}_i \Delta \widetilde{R}_l) = \cup_{i=1}^{i=d} \mathbb{R}_i = \mathbb{R}$ - system components semiring of sets of input signals, which does not include communication functions. Thus the original system of linearly dependent input signals $\widetilde{\mathbb{R}}$, containing a subset of the functions of communication, the subsystem is allocated linearly independent parts of the input signals $\mathbb{R} \subset \widetilde{\mathbb{R}}$, which does not contain subsets containing functions of communication $f_{il}$. Theorem 3 is proved.

Theorem 4. *Let the linear stationary control object contains $d$ linearly dependent inputs, $h$ is linearly dependent O and $r \leq dh$ control channels. On each input and each output simultaneously observed images and archetypes, linear time-invariant convolution type operators, distorted linearly independent multiplicative and additive noise defined on a Hilbert space Besicovitch almost periodic functions. If among all the $d$ preimages running condition $0 \leq |cos(\varphi_{ij})| < 1$ for all inputs $\neq j$, then the frequency control method of noise in systems with independent input actions can be generalized to cases where input actions are linearly dependent signals.*

It is proof. If the input management system observed the initial system of linearly dependent input signals $\widetilde{\mathbb{R}} = \cup_{i=1}^{i=d} \widetilde{\mathbb{R}}_i$, containing a subset of the functions $f_{il}$ connection between inputs $i$ and $l$, then for an arbitrarily selected output $p$ observed process generated by the original linearly dependent system of input actions

$$\widehat{y}_p(t) = y_{xp}(t) + v_p(t) + \beta_p(t) + \vartheta_p(t), \qquad (7)$$

where $y_{xp}(t) = \sum_{c=1}^{c=d} x_c * k_{cp}$ the exact component (image) of the output process generated by the independent components (prototypes) input signals $x_c(t)$, the symbol .*. denotes the operator of convolution type defined in the Hilbert space $B_2(-\infty, +\infty)$.

$v_p(t) = \sum_{r=1}^{r=d-1} \sum_{c=r+1}^{c=d} f_{rc} * (k_{rp} + k_{cp})$ - component of the signal generated at the output $p$ of the functions of communication $f_{rc}$ that occurs between the inputs $r$ and $c$, $k_{rp}$, $k_{cp}$ - the weight functions between the control channels, respectively, inputs $x_c$ $r, c$ and outputs $p$, $(y_p, \bar{v}_p) = 0$ and $x_c$ as and $f_{rc}$ - are linearly independent random processes for all $c$ and $r$,

$\beta_p(t)$ - reduced the additive noise generated by the multiplicative interference, which meets the requirements of independence $(\widehat{y}_p, \bar{\beta}_p) = (v_p, \bar{\beta}_p) = (\beta_p, \bar{\vartheta}_p) = 0$.

$\vartheta_p(t)$ - additive noise at the output of p, satisfying the condition of linear independence ( $\widehat{y}_p, \bar{\vartheta}_p) = (v_p, \bar{\vartheta}_p) = (\beta_p, \bar{\vartheta}_p) = 0$.

Implementation processes at different outlets, for example, $p$ and $r$ are linearly dependent processes $(\widehat{y}_p, \widehat{\bar{y}}_r) \neq 0$ due to mutual coupling between the inputs. However, the components of the process observed at the output $p$ is the sum of mutually orthogonal processes (7). By Theorem 1, the frequency of the harmonic components $\widehat{y}_p(t)$ of the system is a semiring of sets of frequencies

$$\widehat{\Omega}_{yp} = \bigcup_{c=1}^{c=d} \Omega_{xc} \bigcup \Omega_{vp} \bigcup \Omega_{\beta p} \bigcup \Omega_{\vartheta p}. \qquad (8)$$

Suppose that a selection subsystem is linearly independent input signals $\tilde{x}_i = \widehat{\vartheta} x_{xi} + p_{xi} + n_{xi} \in \mathbb{R}_i \subset \mathbb{R}$, each component of which is a finite Fourier series, for example, (3). As a subsystem of $\mathbb{R}$ is a semiring, then by Theorem 1, the system of finite number of sets of frequency harmonic components of the bases $\tilde{x}_i$ is a semicircle, with the representation (6). We find the intersection of the selected linearly independent sets of frequency harmonics, the basic input $i$ (6) and the output $p$ (8). Not difficult to see that $\widetilde{\Omega}_{xi} \cap \widehat{\Omega}_{yp} = \Omega_{xi}$ - set frequency harmonics of the basic components of the exact signal $i$ - m entrance. Knowing the exact set of frequency components of input and output signals $\Omega_{xi}$ by using the Fourier transform of only the frequencies $\omega_k \in \Omega_{xi}$ find accurate projections of the components of input and output signals. Thus, after the separation of independent components of the input signals $\mathbb{R} \subset \widetilde{\mathbb{R}}$ of the original system of linearly dependent input signals $\widetilde{\mathbb{R}}$ frequency control method of noise in systems with independent input actions can be generalized to cases where input actions are linearly dependent signals. Theorem 4 is proved.

7. An example of solving the problem of identification

The proposed method of dealing with noise is used in solving the problems of identification of dynamic characteristics of the multidimensional control objects in a class of linear stationary models. In particular, we solved the problem of identification of dynamic characteristics of the aircraft according to the data obtained in one mode, automatic landing. The task of identifying the dynamic characteristics of the above control channels has been solved by a new method proposed in the following sequence. First by local maxima of the current amplitude spectra [3, pp. 112] observed signals were determined frequency source system sets the frequency of the harmonic components of input signals $\widehat{\Omega}_x = \cup_i \widehat{\Omega}_{xi}$, where $\widehat{\Omega}_{xi} = \widetilde{\Omega}_{xi} \cup (\cup_{k=1}^{k=d-1} (\cup_{m=k+1}^{m=d} \Omega_{fkm}))$, $i = 1,2,3,4$, set of frequencies $\widetilde{\Omega}_{xi}$ determined by the formula (6). $\widehat{\Omega}_y$ - set the frequencies of harmonic components of the output signal $\hat{y}(t)$, represented by the formula (8). Index $p = 1$ is omitted.

Suppose we are given input $i$. There are several ways to define a set of frequencies $\widetilde{\Omega}_{xi}$, which does not contain a subset of the frequencies generated by the communication functions $\Omega_{fik}$, for all inputs $k \neq i$. Consider one of them.

Determined by resolution of the frequency $\delta = 2\pi/T$, where $T$ -duration of implementation. Frequencies are identical if the absolute difference between two frequencies, taken from the sets of frequencies belonging to different inputs, the less resolution. For each selected frequency $\omega_{ik} \in \widehat{\Omega}_{xi}$, $k = 1,2,3,\ldots ci$, where $ci$ - the original dimension of the array of frequencies in the selected input $i$. Are the absolute difference between the frequencies $\Delta_{ikmn} = |\omega_{mn} - \omega_{ik}|$ for all frequencies selected sequentially from the input $m \neq i$. Here $m, n$ - number of inputs, $n, k$ - number of harmonics, $n = 1,2,3,\ldots cm$, where $cm$ - dimension of the original input frequency $m$. If $\Delta_{ikmn} < \delta$, then the system sets the frequency $\widehat{\Omega}_{xi}$ frequency with index $ik$ discarded. In this way we obtain a system of sets of frequencies $\widetilde{\Omega}_{xi}$, which according to the formula (6) does not contain a set of mutually overlapping frequencies generated by the communication functions. The initial dimension of the sets of dependent frequency $\widehat{\Omega}_{x2}$ reduced, for example, the original dimension of the array of frequency dependent components of the inlet throttle position was c2= 137, and the dimension of the array of frequencies $\widetilde{\Omega}_{x2}$ independent harmonic $\tilde{x}_2$ on the same second input has c2= 26.

Next, solve the problem of filtering the frequency harmonics of useful components of the output signals from the frequencies of the harmonics interference. By the corollary of the lemma for a given control channel $\tilde{x}_2 - y$ determined by a set of frequencies of harmonic components of generating and describing the forced motion in the channel controls on the criterion that the frequencies of harmonic components of input and output signals of the selected channel management $\widetilde{\Omega}_{x2y} = \{\omega_\rho = \omega_i \in \widetilde{\Omega}_{x2y}: |\omega_i - \omega_k| \leq \delta, \forall \omega_i \in \widehat{\Omega}_{x2}, \forall \omega_k \in \widehat{\Omega}_y, i \in [1, r_{x2}], k \in [1, r_y]\}$ where $r_y = 35$ - the dimension of the resulting array of harmonic frequencies components of the output signal $\Omega_y, \rho \in [1, q]$, where q - the dimension of the matching frequency harmonic components of generating and describing the forced motion in the channel $\tilde{x}_2 - y$, and in this example, it was found that $q = 6$.

For all the matching frequency $\omega_\rho \in \widetilde{\Omega}_{x2y}$ determined by the Fourier exponents of the input signal $S_{x1}(j\omega_\rho) = a(\omega_\rho) + jb(\omega_\rho)$ and the output signal $S_y(j\omega_\rho) = \gamma(\omega_\rho) + j\beta(\omega_\rho)$. At the lowest frequency $\omega_{\rho 1}$ from the set of frequencies $\omega_{\rho 1} \in \widetilde{\Omega}_{x2y}$ astatizma determines the order in the system for the control action. For what is the value of the frequency transfer function at the lowest frequency $W(j\omega_{\rho 1}) = S_y(j\omega_{\rho 1})/S_x(j\omega_{\rho 1}) = P + jQ$. Depending on the position of the point $W(j\omega_{\rho 1})$ in the complex plane is determined by the order astatizma on the control action. Usually, the specified point can be located in the first, second or third quadrant of the complex plane, which corresponds to the zero, first or second order astatizma $p_a$ [12, s.300].

Obtained by Fourier exponents and order astatizma (in this example $p_a = 1$) can make a number of systems of algebraic equations on the 2 nd to $q$ - th order, which are the Fourier transforms of systems of linear ordinary differential equations with unknown constant coefficients of the 2 nd to $q$ - order. For example, a system of algebraic equations, given the current g$\leq q$ order is represented as

$$\sum_{k=0}^{k=g}(j\omega_{p1})^{k+p_a} T_{k+p_a} S_y(j\omega_{p1}) = S_{x1}(j\omega_{p1}),$$

$$................................$$

$$\sum_{k=0}^{k=g}(j\omega_{pc})^{k+p_a} T_{k+p_a} S_y(j\omega_{pc}) = S_{x1}(j\omega_{pc}),$$

где с$\in [2, g]$.

From this system of equations should be two systems of equations for real and imaginary components. Solving these equations successively from 2 nd to $q$ -th order with respect to coefficients $T_{k+p_a}$, we find that the maximum order of the system of equations $g = q_{max}$ which the system of equations is consistent.

To solve this problem it turned out that the dynamic characteristics of the control channel thrust lever (throttle) - Pitch, taking into account the elastic strain aircraft structure can be described as a first approximation of ordinary differential equations 5-th order with constant coefficients: $T_0 = 0$; $T_1 = -6{,}6537$; $T_2 = -5{,}5133$; $T_3 = -3{,}4828$; $T_4 = -0{,}6220$; $T_5 = -0{,}3525$. There coefficients hew negative signs of the of the differential equation due to the rule of signs taken in the aerodynamics, to the tilt direction throttle and the longitudinal axis of the aircraft in pitch [11, pp. 19]. The differential equation with constant coefficients of the ninth order for the control channel a given pitch - the actual pitch received similarly.

8. Conclusion

Known methods of statistical dynamics is sometimes impossible to apply the solution, for example, the identification problem for the following reasons. Correlation functions cannot be found as well by averaging over time and by averaging over the set of realizations. Stochastic processes are non-stationary processes. They do not possess the ergodic property. The experiments were repeated only under different conditions. Number of inputs exceeds the number of outputs. This eliminates the possibility of obtaining a system of equations. Single realizations of random processes are deterministic functions in such cases.

The concept of measuring a deterministic signal or process is introduced as the sum of one or four linearly independent components. The first component is the exact image or transform operator control channel. The second component is an additive noise. The third component is a reduced additive noise generated by multiplicative noise. This component may be absent in particular cases. The fourth component is a function of communication between the various input signals, if the multidimensional control system is a management system with dependent input actions. All components describe the motion of material bodies on the second law of Newtonian mechanics.

It is assumed that all channels of one-dimensional or multidimensional control systems are described in the first approximation by ordinary differential equations. It follows that the desired function the exact output has continuous derivatives up to $n$ - th order. The exact right-hand side of the differential equation can only be a continuous function. Forced solution exists in the form of a harmonic function on the whole line with harmonic action. Image and transform the differential operator may be provided convergent trigonometric series, and the general solution of differential equation coincides with the forced movement. Such functions are not integrated by Riemann and Lebesgue, and integrated with respect to Stieltjes functions, for example, $g(t) = t/2T$. Function with such properties are well-known elements of the Hilbert space of almost periodic functions of the second order in the sense of Besicovitch In this space defined by a completely continuous operator normal form of convolution type, which is the inverse of differentiation for the exact components of the measured signals. All the components of a signal to be measured belong to the system of sets, which is a semiring. Measure of the semiring generated by the scalar product in Hilbert space of almost periodic functions $B_2$ on the whole line.

In the theory of stochastic processes underlying concept is the concept of uncorrelated random processes. In the theory of nonrandom (deterministic) functions of the basic concept is the notion of a semiring of sets for a system of linearly independent components of the signals. It is shown that if an infinite set of measured signals is a semiring, then the semiring is a finite set of numerical angular frequency harmonic bases in representations of Fourier components of the measured signals. Frequency harmonic basis of the image and transform the operator equations are identical, so the system sets containing the image and inverse image is not a semicircle. These distinguishing features lead to frequency interference mitigation method, an algorithm which is described in the example where, for example, the control system has four inputs and one output.

Initially determined set of frequencies of harmonic components of input and output signals. Determined by the criterion that the frequencies of the resolution. In general, the number of input and output number is determined by the selected control channel. The identification problem is solved sequentially for the selected control channel.

Typically, input signals are not a semicircle. It is shown that the dependent input signal can be a subsystem of independent input signals. Further, the identification problem is solved on a dedicated sub-system of independent input signals. By allocating a subset of matching frequencies for the selected control channel is defined by a subset of the exact frequencies of harmonic components of input and output signals. Fourier exponents are found only on a subset of the frequency of precise components of the signals. So by separating the exact components of the input and output signals of the selected control channel interference. Order and coefficients of the differential equation founds by successively solving systems of equations in the frequency domain for increasing orders, while the condition of compatibility of systems of equations will be satisfied. In this way were found ordinary differential equations of fifth order for the channel lever engine control - pitch plane and the ninth order for the channel - a given - the actual pitch of the aircraft on the data obtained in the regime of one automated aircraft landing, taking into account elastic deformation the airframe.

1.02.2011